# THRESHOLD FOR MONOTONE SYMMETRIC PROPERTIES THROUGH A LOGARITHMIC SOBOLEV INEQUALITY

BY RAPHAËL ROSSIGNOL

*Université René Descartes*

Threshold phenomena are investigated using a general approach, following Talagrand [*Ann. Probab.* **22** (1994) 1576–1587] and Friedgut and Kalai [*Proc. Amer. Math. Soc.* **12** (1999) 1017–1054]. The general upper bound for the threshold width of symmetric monotone properties is improved. This follows from a new lower bound on the maximal influence of a variable on a Boolean function. The method of proof is based on a well-known logarithmic Sobolev inequality on $\{0,1\}^n$. This new bound is shown to be asymptotically optimal.

**1. Introduction.** Threshold phenomena that occur in most discrete probabilistic models have received a lot of attention. One of the archetypal examples is that of the random graphs $\mathcal{G}(n, p(n))$; see [4, 12] or [27]. Consider, for instance, connectivity; see [3]. The probability for $\mathcal{G}(n, p(n))$ to be connected goes from $\varepsilon + o(1)$ to $1 - \varepsilon + o(1)$ when $p(n) = \log n/n + c/n$ and $c$ goes from $\log(1/\log(1/\varepsilon))$ to $\log(1/\log(1/(1-\varepsilon)))$. In this example, the threshold is located around $\log n/n$ and its width is of order $O(1/n)$; see Definition 1.2 below. In the language of statistical physics, threshold phenomena are the "finite-size scaling" parts of phase transitions; see [5]. They have been shown to occur in percolation (see [16]), satisfiability in random constraint models (see, e.g., [5, 10, 14]), local properties in random images (see [9]), reliability (see [21]) and so on. It is therefore of prime interest to find general conditions under which such phenomena occur.

Actually, all the examples cited above can be embedded in the common setting of products of Bernoulli measures on $\{0,1\}^n$; see [15]. Let $n$ be an integer, $p$ a real number in $[0,1]$ and denote by $\mu_{n,p}$ the probability measure on $\{0,1\}^n$ defined by

$$\forall x \in \{0,1\}^n \qquad \mu_{n,p}(x) = p^{\sum_{i=1}^n x_i}(1-p)^{\sum_{i=1}^n(1-x_i)}.$$

We write $\mu_p$ instead of $\mu_{n,p}$ when no confusion is possible.

We are interested in subsets $A$ of $\{0,1\}^n$, the probability $\mu_{n,p}(A)$ of which goes from "almost 0" to "almost 1" over a relatively short interval of values of the probability $p$. The first condition that we shall assume on these subsets is monotonicity.

DEFINITION 1.1. Let $A$ be a subset of $\{0,1\}^n$. The subset $A$ is *monotone* if and only if

$$(x \in A \text{ and } x \preceq y) \implies y \in A,$$

where $\preceq$ is the partial order on $\{0,1\}^n$, defined coordinate-wise.

We shall say that $A$ is *nontrivial* if it is nonempty and different from $\{0,1\}^n$ itself. Let $A$ be a nontrivial monotone subset of $\{0,1\}^n$. It then follows from an elementary coupling technique that the mapping $p \mapsto \mu_p(A)$ is strictly increasing and continuous, thus invertible; see also Lemma 2.2. For $\alpha \in [0,1]$, let $p(\alpha)$ be the unique real number in $[0,1]$ such that $\mu_{p(\alpha)}(A) = \alpha$. The *threshold width* of a subset is the length of the interval over which its probability increases from $\varepsilon$ to $1-\varepsilon$.

DEFINITION 1.2. Let $A$ be a nontrivial monotone subset of $\{0,1\}^n$. Let $\varepsilon \in \,]0, 1/2]$. The threshold width of $A$ at level $\varepsilon$ is

$$\tau(A, \varepsilon) = p(1-\varepsilon) - p(\varepsilon).$$

The first general results on thresholds seem to be those of Margulis [20] and Russo [25], later completed by Talagrand [28, 29]. They related the threshold width to the notion of influence of coordinates. Intuitively, one might say that a subset $A$ will have a narrow threshold unless a few coordinates have a strong influence on its definition [as an example, think of $A = \{x \text{ s.t. } x(1) = 1\}$]. In many cases, this idea is captured by the notion of symmetry.

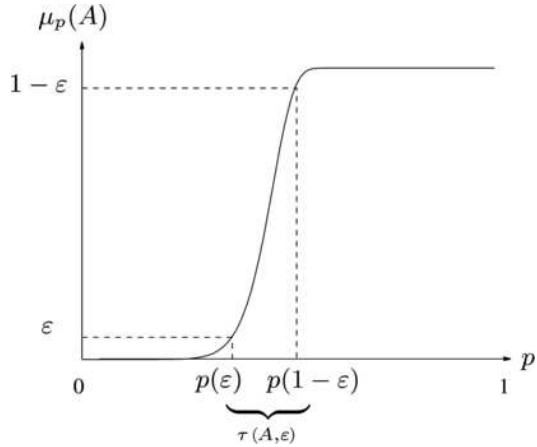

FIG. 1. *Example of a threshold width of level $\varepsilon$.*

DEFINITION 1.3. The subset $A$ of $\{0,1\}^n$ is said to be *symmetric* if and only if there exists a subgroup $G$ of $\mathcal{S}_n$ (group of permutations) acting transitively on $\{1, \ldots, n\}$ such that $A$ is invariant under the action of $G$, that is,

$$\forall g \in G, \ \forall x \in A \qquad g \cdot x = (x_{g^{-1}(1)}, \ldots, x_{g^{-1}(n)}) \in A.$$

This notion of symmetry implies that no coordinate has a stronger influence than any other. It turns out that in most applications, interesting properties are both monotone and symmetric (invariant under permutations of vertices in random graphs, under permutation of clauses in constraint satisfaction problems, etc.). From Corollary 1.4 of [29], one can easily deduce the following theorem that was independently stated by Friedgut and Kalai; see Theorem 3.2 in [15]:

THEOREM 1.4. *There exists a constant $C > 0$ such that, for any nontrivial monotone symmetric subset $A$ of $\{0,1\}^n$ and for all $0 < p < 1$,*

$$\frac{d\mu_p(A)}{dp} \geq C \frac{\log n}{p(1-p)\log(2/(p(1-p)))} \mu_p(A)(1 - \mu_p(A)).$$

It is then easy to derive an upper bound on $\tau(A, \varepsilon)$ from such a result; see Lemma 2.1.

COROLLARY 1.5. *There exists a constant $C > 0$ such that, for any nontrivial monotone symmetric subset $A$ of $\{0,1\}^n$ and for all $0 < \varepsilon < 1/2$,*

$$(1) \qquad \tau(A, \varepsilon) \leq C \sup_{p \in [p(\varepsilon), p(1-\varepsilon)]} \left\{ p(1-p) \log \frac{2}{p(1-p)} \right\} \frac{\log((1-\varepsilon)/\varepsilon)}{\log n}.$$

Corollary 1.5 may in turn be simplified into the following statement:

COROLLARY 1.6. *There exists a constant $C' > 0$ such that, for any nontrivial monotone symmetric subset $A$ of $\{0,1\}^n$ and for all $0 < \varepsilon < 1/2$,*

$$(2) \qquad \tau(A, \varepsilon) \leq C' \frac{\log((1-\varepsilon)/\varepsilon)}{\log n}.$$

Thus, the threshold width of a symmetric monotone property goes to zero as $n$ tends to infinity, and is of order $O(1/\log n)$. When the threshold occurs at a location $p(n)$ which goes to 0 or 1 when $n$ tends to infinity, inequality (2) may be very rough and (1) sharpens this assertion. A natural question regarding these results is whether one can find reasonable bounds for the universal constants $C$ and $C'$.

Both in [29] and [15], the values of $C$ and $C'$ are not explicit. A careful reading of *Talagrand's article* gives the value $C = 120$; see [23], page 23. By following the steps of *Friedgut and Kalai*, the best value that we were able to reach was $C = 5.66$, for a version of Corollary 1.5 where $p(1-p)\log(2/(p(1-p)))$ is replaced by $p(1-p)\log(3/(p(1-p)))$; see [24], page 74. This gives the value $C' = 7.03$. In a recent paper devoted to first passage percolation, [2] gives a new proof of Talagrand's theorem for $p = 1/2$. It is straightforward to generalize this result for any $p \in [0,1]$ and then to deduce a version of Corollary 1.6 with the constant $C' = 3$. Nevertheless, asymptotically, this amounts to twice the best value we offer in this paper.

Our main results are Theorem 4.1 and Corollary 4.3. The first one gives a lower bound on the derivative $d\mu_p(A)/dp$ similar to that of Theorem 4.1, and this bound is asymptotically sharp. Actually, it follows from a slightly more general result on the largest influence of a variable on a Boolean function which we state in Theorem 4.2. Theorem 4.1 implies a sharp version of Corollary 1.5. In particular, we derive a bound for the threshold $\tau(A,\varepsilon)$, similar to that of Corollary 1.6, which is asymptotically equivalent to $(\log((1-\varepsilon)/\varepsilon))/\log n$, thus showing that the universal constant $C'$ can be taken arbitrarily close to 1 for large $n$. These two consequences of Theorem 4.1 are grouped together in Corollary 4.3.

It is tempting to see threshold phenomena as mere consequences of the concentration of product measures, accounted for by a huge variety of probabilistic inequalities; see, for instance, Chapter III in [22], Chapter 2 in [11], the work of Boucheron, Lugosi and Massart [7], Ledoux [19] and Talagrand [30, 31]. Nevertheless, it seems that none of the existing concentration inequalities are able to recover results like Theorem 1.4. The existing proofs of this result all rely on the use of the Beckner–Bonami hypercontractive inequality; see [1, 6]. The main idea of the current article is to replace this central tool by another one, namely a well-known logarithmic Sobolev inequality [see inequality (5)] which allows us to get a sharper result.

Note that another very natural question about the threshold width of a subset $A$ is to what extent it depends on the invariance subgroup $G$ mentioned in Definition 1.3. This question is addressed by Bourgain and Kalai in [8]. Notably, for all but the most basic types of symmetry, the main result of that article asymptotically improves on the bound given in Theorem 4.1. On the other hand, for some "small" symmetry groups (e.g., the cyclic group), Theorem 4.1 is better than the main result in [8].

The paper is organized as follows. Section 2 is devoted to technical results on the derivative of the expectation of a function defined on $\{0,1\}^n$. These results generalize Russo's lemma; see [25] or [16], page 41. The logarithmic Sobolev inequality on which the proof of Theorem 4.1 is based will be explained in Section 3. The proof of the main result is given in Section 4. Finally, the sharpness of Theorem 4.1 is discussed in Section 5.

**2. Threshold width and Russo's lemma.** The usual way to achieve general upper bounds for the threshold width of a set $A$ is to bound $d\mu_p(A)/dp$ below by a suitable function of $p$ and $\mu_p(A)$. To be precise, we will use the following technical lemma:

LEMMA 2.1. *Let $A$ be a monotone, nontrivial subset of $\{0,1\}^n$, $g$ be a continuous positive function on $[0,1]$ and $a$ be a positive real number. The two following propositions are equivalent:*

  (i) $\forall p \in [0,1]$, $\frac{d\mu_p(A)}{dp} \geq \frac{a}{g(p)} \mu_p(A)(1 - \mu_p(A))$;
  (ii) $\forall \alpha \leq \beta \in ]0,1[$, $p_\beta - p_\alpha \leq \frac{1}{a}\sup_{r \in [p(\alpha), p(\beta)]}\{g(r)\} \log \frac{\beta(1-\alpha)}{\alpha(1-\beta)}$.

PROOF. First, let us suppose that (i) is true. Let $\alpha$ and $\beta$ be two real numbers in $]0,1[$ such that $\alpha \leq \beta$. For any $p \in [p(\alpha), p(\beta)]$, we can write

$$\frac{d \log(\mu_p(A)/(1 - \mu_p(A)))}{dp} \geq \frac{a}{\sup_{r \in [p(\alpha); p(\beta)]}\{g(r)\}}.$$

Integrating this inequality between $p(\alpha)$ and $p(\beta)$ the gives (ii). The converse is obtained as follows:

(ii) $\implies$ $\forall \alpha, \beta, 0 < \alpha < \beta < 1$

$$\frac{a}{\sup_{r \in [p(\alpha), p(\beta)]}\{g(r)\}} \leq \frac{\log((\beta(1-\alpha))/(\alpha(1-\beta)))}{p(\beta) - p(\alpha)},$$

$\implies$ $\forall p, q, 0 < p < q < 1$

$$\frac{a}{\sup_{r \in [p,q]}\{g(r)\}} \leq \frac{\log(\mu_q(A)(1 - \mu_p(A))/(\mu_p(A)(1 - \mu_q(A))))}{q - p},$$

which gives (i) by letting $q$ tend to $p$. $\square$

In order to obtain a lower bound for $d\mu_p(A)/dp$, let us define the discrete gradient of a function $f$, from $\{0,1\}^n$ to $\mathbb{R}$:

$$\nabla_i f(x) = f(x_1, \ldots, x_{i-1}, 1, x_{i+1}, \ldots, x_n) - f(x_1, \ldots, x_{i-1}, 0, x_{i+1}, \ldots, x_n).$$

The following lemma is easily obtained by considering the derivative of $\mu_p(x)$ with respect to $p$:

LEMMA 2.2. *For any real function $f$ on $\{0,1\}^n$,*

$$\frac{d}{dp}\int f(x)\,d\mu_p(x) = \sum_{i=1}^{n}\int \nabla_i f(x)\,d\mu_p(x).$$

This expression, when applied to the characteristic function of a monotone set $A$, is equivalent to Russo's lemma; see [16], page 41, or [25]. Indeed, recall the definition of $I_A(i)$, the influence of coordinate $i$ on the subset $A$:

DEFINITION 2.3. Let $n$ be a positive integer and $f$ a function from $\{0,1\}^n$ to $\{0,1\}$. For every $i$ in $\{1,\ldots,n\}$, the influence of variable $i$ on $f$ is the probability of $f$ being nonconstant on the $i$th fiber:

$$I_i(f) = \mu_{n-1,p}(\{x \in \{0,1\}^{n-1},\text{ s.t. } f \text{ is not constant on } l_i(x)\}),$$

where

$$l_i(x) = \{(x_1,\ldots,x_{i-1},u,x_i,\ldots,x_{n-1}) \text{ s.t. } u \in \{0,1\}\}.$$

Let $A$ be a subset of $\{0,1\}^n$. For every $i$ in $\{1,\ldots,n\}$, the *influence of variable $i$ on $A$* is its influence on the characteristic function $\mathbb{1}_A$.

When $f$ is the characteristic function of a monotone set $A$, we have

$$\int \nabla_i f(x)\,d\mu_p(x) = I_A(i).$$

Thus, Lemma 2.2 implies Russo's lemma, which states that for any monotone subset $A$,

$$\frac{d\mu_p(A)}{dp} = \sum_{i=1}^{n} I_A(i).$$

**3. The logarithmic Sobolev inequality on the hypercube.** We introduce (see [29]) the linear operator $\Delta_i$ which acts on any function $f:\{0,1\}^n \to \mathbb{R}$ as follows:

$$\Delta_i f = f - \int f\,d\mu_{1,p}(x_i).$$

This operator is closely related to $\nabla_i$:

(3) $$\Delta_i f(x) = \begin{cases} (1-p)\nabla_i f(x), & \text{if } x_i = 1, \\ -p\nabla_i f(x), & \text{if } x_i = 0. \end{cases}$$

The key property of the operator $\Delta_i$, is that it is the opposite of the generator of a semigroup acting on the $i$th coordinate. To be precise, let us define the semigroup $\{T_t, t \geq 0\}$, acting on $(\{0,1\}, \mu_{1,p})$, of a Markovian jump process with transition rates $p$ from 0 to 1 and $1-p$ from 1 to 0. Its generator $\mathbf{H}$ is the following; see Chapter X in [13]:

$$\mathbf{H}g(x) = \begin{cases} (1-p)(g(0) - g(1)), & \text{if } x = 1, \\ p(g(1) - g(0)), & \text{if } x = 0. \end{cases} \tag{4}$$

Tensorising this semigroup, we obtain a semigroup $\{T_{n,t}, t \geq 0\}$ on $(\{0,1\}^n, \mu_{n,p})$, with generator $\mathbf{L}$:

$$\mathbf{L} = -\sum_{i=1}^{n} \Delta_i.$$

It is known that $\mathbf{H}$ satisfies a logarithmic Sobolev inequality. Let us denote by $\mathbf{Ent}_\mu(f)$ the entropy of a nonnegative function $g$ with respect to a measure $\mu$:

$$\mathbf{Ent}_\mu(g) = \int g \log g \, d\mu - \left(\int g \, d\mu\right) \log\left(\int g \, d\mu\right).$$

The following logarithmic Sobolev inequality, due to Higuchi and Yoshida [18] can be found in [26], Theorem 2.2.8, page 336. For every function $g$ from $\{0,1\}$ to $\mathbb{R}$,

$$\mathbf{Ent}_{\mu_{1,p}}(g) \leq c_{\mathrm{LS}}(p) \int -g\mathbf{H}g \, d\mu_{1,p},$$

where

$$c_{\mathrm{LS}}(p) = \begin{cases} \dfrac{\log(1-p) - \log p}{1 - p - p}, & \text{if } p \neq \dfrac{1}{2}, \\ 2, & \text{if } p = \dfrac{1}{2}. \end{cases}$$

A representation of $p \mapsto c_{\mathrm{LS}}(p)$ is given in Figure 2.

We will now use the tensorization inequality for entropy; see for instance [19]:

$$\mathbf{Ent}_{\mu_{n,p}}(g) \leq \sum_{i=1}^{n} \mathbb{E}_{\mu_{n,p}}(\mathbf{Ent}_{\mu_i}(g)),$$

where $\mathbf{Ent}_{\mu_i}$ means that only the $i$th coordinate is concerned with the integration. This allows us to obtain the following logarithmic Sobolev inequality for any real function $f$ on $\{0,1\}^n$:

$$\mathbf{Ent}_{\mu_{n,p}}(f) \leq c_{\mathrm{LS}}(p) \int -f\mathbf{L}f \, d\mu_{n,p}. \tag{5}$$

In order to see the relevance of inequality (5) in bounding from below the derivative of $p \mapsto \mu_p(A)$, notice now that the term $\int -f\mathbf{L}f \, d\mu_{n,p}$, called the "energy" of the function $f$, is closely related to this derivative if $f = \mathbb{1}_A$. Indeed, whenever $f$ is such that $\nabla_i f \in \{0, 1\}$ for all $i$, Lemma 2.2 can be reformulated as follows:

LEMMA 3.1. *For any function $f$ such that $\nabla_i f \in \{0, 1\}$ for all $i$,*
$$\frac{d \int f \, d\mu_p}{dp} = \frac{1}{p(1-p)} \int -f \mathbf{L} f \, d\mu_{n,p}.$$

PROOF. A simple computation shows how the moments of $\Delta_i f$ and $\nabla_i f$ are related. For any real function $f$ on $\{0, 1\}^n$ and any real number $\alpha \geq 0$,

(6) $$\int |\Delta_i f|^\alpha \, d\mu_p = (p(1-p)^\alpha + (1-p)p^\alpha) \int |\nabla_i f|^\alpha \, d\mu_p.$$

Therefore, as soon as the function $f$ is such that $\nabla_i f \in \{0, 1\}$ for all $i$,
$$\int \nabla_i f(x) \, d\mu_p(x) = \int (\nabla_i f(x))^2 \, d\mu_p(x) = \frac{1}{p(1-p)} \int (\Delta_i f(x))^2 \, d\mu_p(x).$$

This, together with Lemma 2.2, leads to

(7) $$\frac{d}{dp} \int f(x) \, d\mu_p(x) = \frac{1}{p(1-p)} \sum_{i=1}^n \int (\Delta_i f(x))^2 \, d\mu_p(x).$$

Notice that for all functions $f$ and $g$,

(8) $$\int f \Delta_i g \, d\mu_{n,p} = \int \Delta_i f \Delta_i g \, d\mu_{n,p}.$$

Indeed,
$$\int f \Delta_i g \, d\mu_{n,p} - \int \Delta_i f \Delta_i g \, d\mu_{n,p}$$

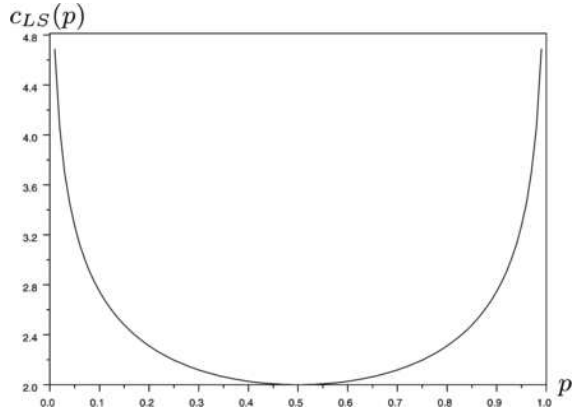

FIG. 2. *The graph of $p \mapsto c_{\mathrm{LS}}(p)$.*

$$= \int \left( \int f \, d\mu_{1,p}(x_i) \right) \left( g - \int g \, d\mu_{1,p}(x_i) \right) d\mu_{n,p}$$

$$= \int \left( \int f \, d\mu_{1,p}(x_i) \right) \left( \int \left( g - \int g \, d\mu_{1,p}(x_i) \right) d\mu_{1,p}(x_i) \right) d\mu_{n,p}$$

$$= 0.$$

Therefore, from equation (7) we obtain

$$\frac{d}{dp} \int f(x) \, d\mu_p(x) = \frac{1}{p(1-p)} \sum_{i=1}^n \int f(x) \Delta_i f(x) \, d\mu_p(x)$$

$$= \frac{1}{p(1-p)} \int f(x) \sum_{i=1}^n \Delta_i f(x) \, d\mu_p(x),$$

which leads to the desired result. $\square$

The role played by the logarithmic Sobolev inequality (5) in the subsequent proof is very similar to the one played by hypercontractivity of the same semigroup in the result of Talagrand [29]. In that article, hypercontractivity for the semigroup $\{T_{n,t}, t \geq 0\}$ is achieved from the $p = 1/2$ case by using a symmetrization technique. Actually, a theorem due to Gross [17] gives an exact equivalence between hypercontractivity and the existence of a logarithmic Sobolev inequality. But the hypercontractivity function found by Talagrand is not optimal. Indeed, the one obtained by using Gross' theorem and inequality (5) is better (and optimal). Notice, though, that when Talagrand's article was published in 1994, the precise logarithmic Sobolev constant $c_{\mathrm{LS}}(p)$ was not yet known.

We finish this section by recalling a classical Poincaré inequality on $\{0,1\}^n$ that will be useful in the sequel. Let $g$ be a function on $\{0,1\}$. A simple computation relates the variance of $g$ and the energy of $g$ associated to **H**:

$$\mathrm{Var}_{\mu_{1,p}}(g) = \int -g\mathbf{H}g \, d\mu_{1,p}.$$

The Jensen inequality implies the following tensorization property for the variance of a function $f$ from $\{0,1\}^n$ to $\mathbb{R}$; see [19]:

$$\mathrm{Var}_{\mu_{n,p}}(f) \leq \sum_{i=1}^n \mathbb{E}_{\mu_{n,p}}(\mathrm{Var}_{\mu_i}(f)),$$

where $\mathrm{Var}_{\mu_i}$ means that only the $i$th coordinate is concerned with the integration. This leads to the following Poincaré inequality for any real function $f$ on $\{0,1\}^n$:

(9) $$\mathrm{Var}_{\mu_{n,p}}(f) \leq \int -f\mathbf{L}f \, d\mu_{n,p}.$$

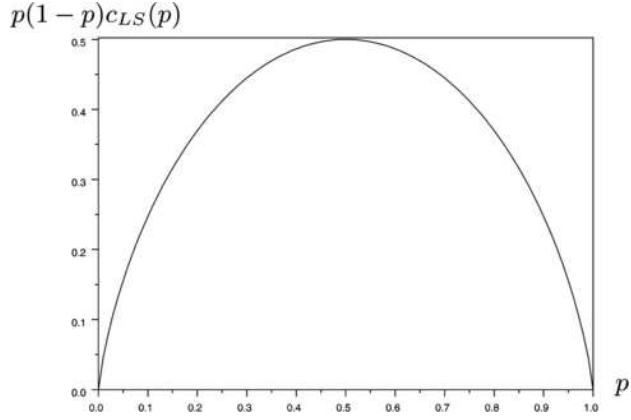

FIG. 3.  *The graph of $p \mapsto p(1-p)c_{\mathrm{LS}}(p)$.*

**4. Main result.** We now turn to the statement of Theorem 4.1, the main result of this article.

THEOREM 4.1.  *Let $s(n)$ be the following sequence of real numbers:*

$$s(n) = \log n - \max\left\{\log\left(\frac{e^{2+4/e}}{2^{5+4/e}}\left(\log \frac{n}{(\log n)^2}\right)^{3+4/e}\right), 2\log(\log n)\right\}.$$

*For every integer $n \geq 2$, every real number $p \in ]0,1[$ and every nontrivial monotone symmetric subset $A$ of $\{0,1\}^n$,*

$$\forall p \in ]0,1[ \qquad \frac{d\mu_p(A)}{dp} \geq \frac{s(n)}{p(1-p)c_{\mathrm{LS}}(p)}\mu_p(A)(1-\mu_p(A)).$$

A graph of $p \mapsto p(1-p)c_{\mathrm{LS}}(p)$ is shown in Figure 3. Also, one can check numerically that

$$\forall n \geq 2 \qquad s(n) > 0$$

and

$$\forall n \geq 275 \qquad \log\left(\frac{e^{2+4/e}}{2^{5+4/e}}\left(\log \frac{n}{(\log n)^2}\right)^{3+4/e}\right) \geq 2\log(\log n).$$

Therefore,

$$\forall n \geq 275 \qquad s(n) \geq \log n - \left(3 + \frac{4}{e}\right)\log\log\frac{n}{(\log n)^2}.$$

Of course, as $n$ tends to infinity, $s(n)$ is equivalent to $\log n$.

Actually, Theorem 4.1 is an easy consequence of the following, slightly more general, result on the largest influence of a variable on a Boolean function (see Definition 2.3):

THEOREM 4.2. *For every integer $n \geq 2$, every real number $p \in \,]0,1[$ and every function $f$ from $\{0,1\}^n$ to $\{0,1\}$, the largest influence of a variable on $f$ is bounded below as follows:*

$$\max\{I_j(f) \ s.t. \ j=1,\ldots,n\} \geq \frac{\mathrm{Var}(f)s(n)}{np(1-p)c_{\mathrm{LS}}(p)}.$$

PROOFS OF THEOREMS 4.1 AND 4.2. Let $f$ be a function on $\{0,1\}^n$ with values in $\mathbb{R}$ and define

$$\forall j \in \{1,\ldots,n\} \qquad V_j = \mathbb{E}[f|(x_1,\ldots,x_j)] - \mathbb{E}[f|(x_1,\ldots,x_{j-1})].$$

Then, we write $f - \mathbb{E}(f)$ as a sum of the martingale increments $V_j$, for $j = 1, \ldots, n$:

$$f - \mathbb{E}[f] = \sum_{j=1}^{n} V_j.$$

Notice that martingale increments are always orthogonal:

$$\forall j \neq k \qquad \int V_j V_k \, d\mu_p = 0.$$

In addition, since $\int V_j \, dx_i$ and $\int V_k \, dx_i$ are two different martingale increments for another filtration, they are also orthogonal:

$$\forall i, \ \forall j \neq k \qquad \int \left( \int V_j \, dx_i \int V_k \, dx_i \right) d\mu_p = 0.$$

We apply the logarithmic Sobolev inequality (5) to each increment $V_j$:

$$c_{\mathrm{LS}}(p) \int V_j \sum_{i=1}^{n} \Delta_i V_j \, d\mu_p \geq \mathbf{Ent}_{\mu_p}(V_j^2)$$

$$= \int V_j^2 \log V_j^2 \, d\mu_p - \int V_j^2 \, d\mu_p \log \int V_j^2 \, d\mu_p.$$

Summing these inequalities for $j = 1, \ldots, n$ results in the following:

$$c_{\mathrm{LS}}(p) \sum_{i=1}^{n} \sum_{j=1}^{n} \int V_j \Delta_i V_j \, d\mu_p$$

$$\geq \sum_{j=1}^{n} \int V_j^2 \log V_j^2 \, d\mu_p + \sum_{j=1}^{n} \|V_j\|_2^2 \log \frac{1}{\|V_j\|_2^2}.$$

We now claim that the sum of the energies of the increments $V_j$ is equal to the energy of $f$:

(10) $$\sum_{j=1}^{n} \int V_j \Delta_i V_j \, d\mu_p = \int f \Delta_i f \, d\mu_p.$$

Indeed,

$$\int f \Delta_i f \, d\mu_p = \int \sum_{j=1}^{n} V_j \Delta_i \sum_{k=1}^{n} V_k \, d\mu_p$$

$$= \sum_{j=1}^{n} \int V_j \Delta_i V_j \, d\mu_p + \sum_{j \neq k} \int V_j \Delta_i V_k \, d\mu_p.$$

Recall that

$$\Delta_i V_k = V_k - \int V_k \, dx_i.$$

Thus,

$$\int V_j \Delta_i V_k \, d\mu_p = \int V_j V_k \, d\mu_p - \int V_j \left( \int V_k \, dx_i \right) d\mu_p$$

$$= \int V_j V_k \, d\mu_p - \int \left( \int V_j \, dx_i \right) \left( \int V_k \, dx_i \right) d\mu_p,$$

and each term of the last sum is null whenever $j \neq k$. This proves the claim (10).

One can now write

$$c_{\mathrm{LS}}(p) \sum_{i=1}^{n} \int f \Delta_i f \, d\mu_p \geq \sum_{j=1}^{n} \int V_j^2 \log V_j^2 \, d\mu_p + \sum_{j=1}^{n} \|V_j\|_2^2 \log \frac{1}{\|V_j\|_2^2}.$$

From equation (8), we deduce

$$\int f \Delta_i f \, d\mu_p = \int (\Delta_i f)^2 \, d\mu_p,$$

and therefore,

$$c_{\mathrm{LS}}(p) \sum_{i=1}^{n} \int (\Delta_i f)^2 \, d\mu_p \geq \underbrace{\sum_{j=1}^{n} \int V_j^2 \log V_j^2 \, d\mu_p}_{(1)} + \underbrace{\sum_{j=1}^{n} \|V_j\|_2^2 \log \frac{1}{\|V_j\|_2^2}}_{(2)}.$$

First, let us rewrite $V_j$ as follows:

(11) $$V_j = \mathbb{E}[-\Delta_j f | x_1, \ldots, x_j].$$

Using Jensen's inequality,

$$\|V_j\|_2^2 \leq \|\Delta_j f\|_2^2.$$

Let us note that

$$\Delta = \max_{j} \|\Delta_j f\|_2^2.$$

We can then obtain a lower bound for the term (2) as follows:

$$(2) = \sum_{j=1}^{n} \|V_j\|_2^2 \log \frac{1}{\|V_j\|_2^2} \geq \sum_{j=1}^{n} \|V_j\|_2^2 \log \frac{1}{\|\Delta_j f\|_2^2}$$

$$\geq \mathrm{Var}(f) \log \frac{1}{\Delta}. \tag{12}$$

Let us split each term of the sum (1) in the following way:

$$\int V_j^2 \log V_j^2 \, d\mu_p = \underbrace{\int V_j^2 \log V_j^2 \mathbb{1}_{V_j^2 \leq t} \, d\mu_p}_{(1\mathrm{a})} + \underbrace{\int V_j^2 \log V_j^2 \mathbb{1}_{V_j^2 > t} \, d\mu_p}_{(1\mathrm{b})}.$$

Since the function $x \mapsto x \log x$ is nonincreasing on $[0, 1/e]$, we can write, for every $t \leq 1/e^2$,

$$(1\mathrm{a}) = \int 2|V_j| \times |V_j| \log |V_j| \mathbb{1}_{V_j^2 \leq t} \, d\mu_p \geq \sqrt{t} \log t \int |V_j| \mathbb{1}_{V_j^2 \leq t} \, d\mu_p$$

$$\geq \sqrt{t} \log t \int |V_j| \, d\mu_p,$$

since $\sqrt{t} \log t$ is nonpositive.

From equation (11), and using Jensen's inequality, we derive the following:

$$\int |V_j| \, d\mu_p \leq \int |\Delta_j f| \, d\mu_p.$$

We now use equation (6) with $\alpha = 1$, then the fact that $\nabla_j f \in \{0, 1\}$ and finally equation (6) with $\alpha = 2$:

$$\int |\Delta_j| \, d\mu_p = 2p(1-p) \int |\nabla_j f| \, d\mu_p$$

$$= 2p(1-p) \int (\nabla_j f)^2 \, d\mu_p$$

$$= 2 \int (\Delta_j f)^2 \, d\mu_p.$$

Moreover, the log function being increasing, we have

$$(1\mathrm{b}) \geq \log t \int V_j^2 \, d\mu_p.$$

Summing the lower bounds thus collected, we find

$$(1) \geq 2\sqrt{t} \log t \sum_{j=1}^{n} \int -f \mathbf{L}_j f \, d\mu_p + \log(t) \sum_{j=1}^{n} \int V_j^2 \, d\mu_p$$

$$= 2I\sqrt{t} \log t + \mathrm{Var}(f) \log(t), \tag{13}$$

where we have introduced the notation

$$I = \sum_{i=1}^{n} \int (\Delta_i f)^2 \, d\mu_p.$$

We would like to choose $t$ so as to maximize expression (13). It is easier to equalize the terms $2I\sqrt{t}\log t$ and $\text{Var}(f)\log(t)$. We would therefore be tempted to take

$$t = \left(\frac{\text{Var}(f)}{2I}\right)^2,$$

but we have to maintain agreement with the hypothesis that $t \leq 1/e^2$, whereas we only know, by the Poincaré inequality (9), that

$$I \geq \text{Var}(f).$$

Let us choose, then,

$$t = \left(\frac{\text{Var}(f)}{eI}\right)^2.$$

Thus,

(14) $$(1) \geq \text{Var}(f)\log\left(\frac{\text{Var}(f)}{eI}\right)^{2+4/e}.$$

Collecting lower bounds on (1) and (2) from (12) and (14) and using the trivial bound, we have

$$\Delta \geq \frac{I}{n},$$

so we get

(15) $$c_{\text{LS}}(p)I \geq \text{Var}(f)\log\left(\left(\frac{\text{Var}(f)}{eI}\right)^{2+4/e}\frac{1}{\Delta}\right),$$

$$c_{\text{LS}}(p)\Delta \geq \frac{1}{n}\text{Var}(f)\log\left(\left(\frac{\text{Var}(f)}{eI}\right)^{2+4/e}\frac{1}{\Delta}\right).$$

Now, let us consider the following disjunction:
• Either

$$c_{\text{LS}}(p)I \geq \text{Var}(f)\log\frac{n}{(\log n)^2},$$

and therefore,

$$c_{\text{LS}}(p)\Delta \geq \frac{1}{n}\text{Var}(f)\log\frac{n}{(\log n)^2},$$

- or

$$c_{\mathrm{LS}}(p)I < \mathrm{Var}(f)\log\frac{n}{(\log n)^2},$$

and thus, using (15),

(16) $$c_{\mathrm{LS}}(p)\Delta \geq \frac{1}{n}\mathrm{Var}(f)\log\left(\left(\frac{c_{\mathrm{LS}}(p)}{e\log(n/(\log n)^2)}\right)^{2+4/e}\frac{1}{\Delta}\right).$$

Then, again, we either have

$$c_{\mathrm{LS}}(p)\Delta \geq \frac{1}{n}\mathrm{Var}(f)\log\frac{n}{(\log n)^2},$$

or,

$$c_{\mathrm{LS}}(p)\Delta < \frac{1}{n}\mathrm{Var}(f)\log\frac{n}{(\log n)^2},$$

which gives, via inequality (16),

$$c_{\mathrm{LS}}(p)\Delta \geq \frac{1}{n}\mathrm{Var}(f)\log\left(\left(\frac{c_{\mathrm{LS}}(p)}{e\log(n/(\log n)^2)}\right)^{2+4/e}\frac{c_{\mathrm{LS}}(p)}{\mathrm{Var}(f)\log(n/(\log n)^2)}\right).$$

In any case,

(17) $$\Delta \geq \frac{\mathrm{Var}(f)}{nc_{\mathrm{LS}}(p)}$$
$$\times \min\left\{\log\frac{n}{(\log n)^2}, \log\left(\frac{nc_{\mathrm{LS}}(p)^{3+4/e}}{e^{2+4/e}\mathrm{Var}(f)(\log(n/(\log n)^2))^{3+4/e}}\right)\right\}.$$

Notice now that

$$\Delta = p(1-p)\max\{I_j(f) \text{ s.t. } j=1,\ldots,n\}$$

and, of course, when $f$ is a Boolean function on $n$ variables,

$$\mathrm{Var}(f) \leq \tfrac{1}{4}.$$

Therefore, inequality (17), together with the observation that $c_{\mathrm{LS}}(p) \geq 2$, leads to

$$\max\{I_j(f) \text{ s.t. } j=1,\ldots,n\} \geq \frac{\mathrm{Var}(f)s(n)}{np(1-p)c_{\mathrm{LS}}(p)}.$$

The proof of Theorem 4.2 is complete. To see how this implies Theorem 4.1, let $f = \mathbb{1}_A$ be the characteristic function of a monotone symmetric subset $A$ of $\{0,1\}^n$. Since $f$ is a symmetric function, the influences of $f$ are all equal and thus,

$$\max\{I_j(f) \text{ s.t. } j=1,\ldots,n\} = \frac{1}{n}\sum_{j=1}^n I_j(f)$$
$$= \frac{1}{n}\frac{d\mu_p(A)}{dp},$$

where the last equality follows from Lemma 2.2 and the fact that $A$ is monotone. Therefore, Theorem 4.2 applied to the Boolean function $f$ implies that

$$\frac{d\mu_p(A)}{dp} \geq \frac{\text{Var}(f)s(n)}{p(1-p)c_{\text{LS}}(p)}. \qquad \square$$

We now turn to the upper bound on the threshold width of a nontrivial symmetric set.

COROLLARY 4.3. *Let $s(n)$ be defined as in Theorem 4.1. For every integer $n \geq 2$, every real number $\varepsilon \in {]}0,1/2{[}$ and every nontrivial monotone symmetric subset $A$ of $\{0,1\}^n$,*

$$(18) \qquad \tau(A,\varepsilon) \leq \sup_{p \in [p(\varepsilon), p(1-\varepsilon)]} \{p(1-p)c_{\text{LS}}(p)\}\frac{\log((1-\varepsilon)/\varepsilon)}{2s(n)}$$

*and, in particular,*

$$(19) \qquad \tau(A,\varepsilon) \leq \frac{\log((1-\varepsilon)/\varepsilon)}{s(n)}.$$

PROOF. Theorem 4.1 ensures that

$$\frac{d\mu_p(A)}{dp} \geq \frac{\text{Var}(f)s(n)}{p(1-p)c_{\text{LS}}(p)}.$$

Since $s(n)$ is positive for all $n \geq 2$, inequality (18) follows from Lemma 2.1. Notice that $p(1-p)c_{\text{LS}}(p) \leq 1/2$ (see Figure 3). This implies inequality (19). $\square$

Recalling that as $n$ tends to infinity, $s(n)$ is equivalent to $\log n$, the second assertion of Corollary 4.3 means that, asymptotically, we can lower by a factor of 2 the best constant in Friedgut and Kalai's theorem ([15], Corollary 1.6), obtained by following the work of Benjamini, Kalai and Schramm [2].

**5. Sharpness of the bound.** Let us discuss now the sharpness of Theorem 4.1 and its corollaries. The following lemma implies that Theorem 4.1 is optimal, if the desired lower bound involves $\mu_p(A)(1-\mu_p(A))$ and $p(1-p)c_{\mathrm{LS}}(p)$, or some equivalents, as $\mu_p(A)$ or $p$ tends to zero:

LEMMA 5.1. *Suppose that there exist two positive functions $f$ and $g$ and a sequence of positive real numbers $\{a(n), n \in \mathbb{N}^*\}$ such that for every $n \in \mathbb{N}^*$, every monotone symmetric subset $A \subset \{0,1\}^n$ and every $p \in ]0,1[$,*

$$\frac{d\mu_p(A)}{dp} \geq \frac{a(n)}{g(p)} f(\mu_p(A)),$$

*with*

$$\frac{f(x)}{x} \xrightarrow{x \to 0} 1 \quad and \quad \frac{g(p)}{p \log(1/p)} \xrightarrow{p \to 0} 1.$$

*Then,*

$$\limsup_{n \to +\infty} \frac{a(n)}{\log n} \leq 1.$$

PROOF. For $n \geq 2$, consider the following monotone symmetric subset:

$$B_n = \{x \in \{0,1\}^n \text{ s.t. } \exists i \in \{1,\ldots,n\}, x_i = 1\}.$$

The probability of $B_n$ is

$$\mu_p(B_n) = 1 - (1-p)^n.$$

Therefore,

$$\frac{d\mu_p(B_n)}{dp} = n(1-p)^{n-1}.$$

Fix $\varepsilon \in ]0, 1/2[$. Suppose now that $p = p(n)$ is such that $\mu_p(B_n) = \varepsilon$. Then $p(n)$ tends to zero as $n$ tends to infinity. Therefore,

$$(1-p(n))^{n-1} = 1 - \varepsilon + o(1),$$

$$np(n) = \log \frac{1}{1-\varepsilon} + o(1).$$

Thus,

$$\left.\frac{d\mu_p(B_n)}{dp}\right|_{p=p(n)} = \frac{1}{p(n)}(1-\varepsilon)\log\frac{1}{1-\varepsilon} + o\left(\frac{1}{p(n)}\right)$$

and

$$\log\left(\frac{1}{p(n)}\right) = \log \frac{n}{\log(1/(1-\varepsilon))} + o(1).$$

Hence,

$$\left.\frac{d\mu_p(B_n)}{dp}\right|_{p=p(n)} = \frac{\log n}{p(n)\log(1/p(n))}(1-\varepsilon)\log\frac{1}{1-\varepsilon} + o\left(\frac{\log n}{p(n)\log(1/p(n))}\right).$$

Therefore,

$$\lim_{n\to+\infty} \left.\frac{d\mu_p(B_n)}{dp}\right|_{p=p(n)} \times \frac{1}{\log n}\frac{p(n)\log(1/p(n))}{(1-\varepsilon)\log(1/(1-\varepsilon))} = 1.$$

Suppose now that there exist two positive functions $f$ and $g$ and a sequence of positive real numbers $\{a(n), n\in\mathbb{N}\}$ such that for every $n\in\mathbb{N}$, every monotone symmetric subset $A\subset\{0,1\}^n$ and every $p\in\,]0,1[$, we have

(20) $$\frac{d\mu_p(A)}{dp} \geq \frac{a(n)}{g(p)}f(\mu_p(A))$$

and

$$\frac{f(x)}{x} \xrightarrow{x\to 0} 1 \quad \text{and} \quad \frac{g(p)}{p\log(1/p)} \xrightarrow{p\to 0} 1.$$

Inequality (20) holds, in particular, for $A=B_n$ and $p=p(n)$. Therefore,

$$\begin{aligned}
1 &= \lim_{n\to+\infty} \left.\frac{d\mu_p(B_n)}{dp}\right|_{p=p(n)} \times \frac{1}{\log n}\frac{p(n)\log(1/p(n))}{(1-\varepsilon)\log(1/(1-\varepsilon))} \\
&\geq \limsup_{n\to+\infty} \frac{a(n)}{\log n} \times \frac{p(n)\log(1/p(n))}{g(p(n))} \times \frac{f(\varepsilon)}{(1-\varepsilon)\log(1/(1-\varepsilon))} \\
&= \frac{f(\varepsilon)}{(1-\varepsilon)\log(1/(1-\varepsilon))} \limsup_{n\to+\infty} \frac{a(n)}{\log n}.
\end{aligned}$$

This inequality is valid for any $\varepsilon\in\,]0,1[$. Since $\frac{f(\varepsilon)}{(1-\varepsilon)\log 1/(1-\varepsilon)}$ tends to one as $\varepsilon$ goes to zero,

$$\limsup_{n\to+\infty} \frac{a(n)}{\log n} \leq 1. \qquad\square$$

As suggested by Lemma 2.1, one can see that inequality (18) in Corollary 4.3 is also asymptotically sharp. Nevertheless, it remains unknown whether inequality (19) is optimal or not. Indeed, this inequality is equivalent to equality (18) only when the threshold is located at $p=1/2$. Following [15] in studying the "Tribes example," it is possible to construct a sequence of monotone symmetric subsets $C_n\subset\{0,1\}^n$ with a threshold located in $1/2$ and such that, for all $\varepsilon$ in $]0,1/2[$,

$$\tau(C_n,\varepsilon) = \frac{\log 2(\log\log(1/(1-\varepsilon)) - \log\log(1/\varepsilon))}{\log n} + o\left(\frac{1}{\log n}\right).$$

When $\varepsilon$ tends to $1/2$, this threshold width gets close to $\log((1-\varepsilon)/\varepsilon)/2\log n$. Therefore, it remains an open problem to find an optimal upper bound for the threshold width of a symmetric property whose threshold is located at $1/2$.

**Acknowledgments.** I would like to thank Bernard Ycart for introducing me to the threshold phenomena and for his patient guidance throughout my dealing with this topic. I am grateful to an anonymous referee for many helpful comments and notably for suggesting the more general result on influences underlying a first version of Theorem 4.1. This lead to the statement of Theorem 4.2.

LABORATOIRE MAP5
UNIVERSITÉ RENÉ DESCARTES–PARIS 5
45 RUE DES SAINTS-PÈRES
75270 PARIS CEDEX 06
FRANCE
E-MAIL: rost@math-info.univ-paris5.fr
URL: http://www.math-info.univ-paris5.fr/~rost/